\documentclass[a4paper]{amsart}
\usepackage[utf8]{inputenc}
\usepackage[T1]{fontenc}
\pdfoutput=1

\usepackage[style=alphabetic,maxbibnames=10,backend=biber,isbn=false,giveninits,date=year]{biblatex}
\DeclareSourcemap{
	\maps{
		\map{
			\step[fieldset=pagetotal, null]
			\step[fieldset=address, null]
		}
	}
}
\appto{\bibsetup}{\emergencystretch=0.75em}
\usepackage{paralist}
\usepackage{hyperref}

\makeatletter
\newcommand{\proofstep}[1]{%
  \par%
  \addvspace{\medskipamount}%
  \textit{#1\@addpunct{.}}\enspace\ignorespaces
}  
  
\makeatother

  {\par%
   \hspace{\parindent}\textsc{Acknowledgements.}\enskip\ignorespaces}
  {\par}
  
\newcommand{\doublerightarrow}{%
  \mathrel{%
  \mathop{\vcenter{%
    \offinterlineskip\ialign{##\cr
    $\rightarrow$\cr\noalign{\kern.5ex}
    $\rightarrow$\cr}}}}
}

\usepackage{amsmath, amssymb, amsfonts, mathtools, amsthm, mathrsfs, dsfont, fge, extpfeil}
\DeclareFontFamily{U}{matha}{}
\DeclareFontShape{U}{matha}{m}{n}{
  <-5.5>    matha5
  <5.5-6.5> matha6 
  <6.5-7.5> matha7
  <7.5-8.5> matha8
  <8.5-9.5> matha9
  <9.5-11>  matha10
  <11->     matha12
}{}

\usepackage{pifont}
\makeatletter
\newcommand\Pimathsymbol[3][\mathord]{%
  #1{\@Pimathsymbol{#2}{#3}}}
\def\@Pimathsymbol#1#2{\mathchoice
  {\@Pim@thsymbol{#1}{#2}\tf@size}
  {\@Pim@thsymbol{#1}{#2}\tf@size}
  {\@Pim@thsymbol{#1}{#2}\sf@size}
  {\@Pim@thsymbol{#1}{#2}\ssf@size}}
\def\@Pim@thsymbol#1#2#3{%
  \mbox{\fontsize{#3}{#3}\Pisymbol{#1}{#2}}}
\makeatother

\DeclareFontFamily{U}{mathb}{\hyphenchar\font45}
\DeclareFontShape{U}{mathb}{m}{n}{
  <-6> mathb5 <6-7> mathb6 <7-8> mathb7
  <8-9> mathb8 <9-10> mathb9
  <10-12> mathb10 <12-> mathb12
  }{}

\usepackage{tikz,tikz-cd,multirow}
\usetikzlibrary{matrix,arrows,positioning,calc,fit}
\usepackage{resizegather}

\newcommand\eps{\varepsilon}
\newcommand\B{\mathcal{B}}
\newcommand\M{\mathcal{M}}
\newcommand\N{\mathbb{N}}
\newcommand\cN{\mathcal{N}}
\newcommand\Sp{\mathrm{Sp}}
\newcommand\Ord{Q_1^\mathrm{ord}}

\DeclareMathOperator{\embeds}{\hookrightarrow}
\DeclarePairedDelimiter{\set}{\{}{\}}
\DeclarePairedDelimiter{\card}{|}{|}
\def\abs#1{\left\lvert #1 \right\rvert}
\DeclarePairedDelimiter{\spann}{\langle}{\rangle}

\DeclareMathOperator{\im}{im}
\DeclareMathOperator{\mods}{mod}
\DeclareMathOperator{\ind}{ind}
\DeclareMathOperator{\dir}{dir}
\newcommand{\edge}{\mathrel{\relbar\mkern-9mu\relbar}}
\newcommand{\from}{\leftarrow}
\newcommand{\rightloop}{\Pimathsymbol[\mathrel]{mathb}{"FD}}
\newcommand\isom{\mathrel{\cong}}
\def\suchthat{\mathrel{|}}
\def\setminus{\mathbin{\fgebackslash}}

\def\defined#1{\textbf{#1}}

\theoremstyle{definition}
\newtheorem{mydef}{Definition}[section]

\theoremstyle{plain}
\newtheorem{theorem}[mydef]{Theorem}
\newtheorem*{theorem*}{Theorem}

\newtheorem{proposition}[mydef]{Proposition}
\newtheorem{lemma}[mydef]{Lemma}
\newtheorem{corollary}[mydef]{Corollary}
\newtheorem*{corollary*}{Corollary}

\theoremstyle{remark}
\newtheorem{example}[mydef]{Example}

\newtheorem*{example*}{Example}
\newtheorem*{examples*}{Examples}

\title{Clannish algebras are of amenable representation type}
\author{Sebastian Eckert}
\address{Sebastian Eckert, Max-Planck-Institut für Mathematik, Vivatsgasse 7, 53111~Bonn, Germany.}
\email{eckert@mpim-bonn.mpg.de}

\subjclass[2020]{16G20,16G60,05C10}
\keywords{representations of finite dimensional algebras, amenable representation type, string algebras, clannish algebras, coefficient quivers, fragmentable graphs}

\begin{document}

\begin{abstract}
We revisit G. Elek's notion of amenable representation type, where algebras are characterised by every indecomposable module being ``almost'' the direct sum of modules of bounded dimension. We give a new proof of his result that string algebras are amenable that relies on the planarity of the coefficient quivers of indecomposable string and band modules. Using this connection to graph theory, we then show that clannish algebras are also of amenable type.
\end{abstract}

\maketitle

\section{Introduction}

\begin{mydef} \label{def:hyperfinite} \label{def:amenability}
Let $k$ be a field, $A$ a finite dimensional $k$-algebra and $\M \subseteq \mods A$ a family of $A$-modules. We say that $\M$ is \defined{hyperfinite}\index{hyperfinite} provided for every $\varepsilon > 0$ there exists $L_\varepsilon > 0$ such that for every $M \in \M$ there exists a submodule $N \subseteq M$ such that
\begin{equation} \label{eq:HFSubmoduleBig} \dim_k N \geq (1-\varepsilon) \dim_k M, \end{equation}
and there are modules 
\begin{equation} N_1, N_2, \dots, N_t \in \mods A, \, \text{with} \, \dim_k N_i \leq L_\varepsilon,\end{equation}
such that 
\begin{equation} N \isom \bigoplus_{i=1}^{t} N_i.\end{equation}

A $k$-algebra $A$ is said to be of \defined{amenable representation type}\index{representation type!amenable} provided $\mods A$ itself is a hyperfinite family.
\end{mydef}

\begin{proposition} \label{prop:ExtendingHFfromSubmodulesOfBoundedCodimension}
Let $A$ be a finite dimensional $k$-algebra and $\M,\cN \subseteq \mods A$ where $\cN$ is hyperfinite. If there is some $H \geq 0$ such that for all $M \in \M$, there exists a submodule $N \subseteq M$ with $N \in \cN$, of codimension less than or equal to $H$, then $\M$ is also hyperfinite.
\begin{proof}
Let $H \geq 0$ and $\varepsilon > 0$. If the dimension of the modules in $\M$ was bounded, say by $K$, we can set $L_\varepsilon \coloneqq K$ and choose $N = M$  for all $M \in \M$, and we are done.
On the other hand, if the dimension is not bounded, there is $M \in \M$ with 
$\dim_k M > \frac{2H}{\varepsilon}$. We choose a submodule $N \in \cN$ of codimension bounded by $H$. Since $\cN$ is hyperfinite, there is some submodule $Y \subseteq N$ such that $\dim Y \geq (1-\frac{\varepsilon}{2})\dim N$, while $Y$ decomposes into direct summands of dimension less than or equal to $L^{\cN}_{\frac{\varepsilon}{2}}$.

We thus have that
\begin{align*}
\dim Y &\geq \left(1-\frac{\varepsilon}{2}\right) \dim N = \dim N - \frac{\varepsilon}{2} \dim N\\
&\geq \left( \dim M - H \right) - \frac{\varepsilon}{2} \dim M \\
&\geq \dim M - \frac{\varepsilon}{2} \dim M - \frac{\varepsilon}{2} \dim M\\
&= (1-\varepsilon) \dim M,
\end{align*}
using that we have $H \leq \frac{\varepsilon}{2} \dim M$.
What is more, $Y$ decomposes into direct summands of dimension less than or equal to $L^{\cN}_{\frac{\varepsilon}{2}}$. If we therefore choose $L^{\M}_\varepsilon$ to be the maximum of $L^{\cN}_{\frac{\varepsilon}{2}}$ and $\frac{2H}{\varepsilon}$, we have shown that $\M$ is hyperfinite.
\end{proof}
\end{proposition}

\section{Planarity of coefficient quivers}

\subsection{Coefficient quivers and tree modules}

In the following, we will consider the path algebra of a given quiver $Q$. 
Recall from \cite[Section~1]{Ringel1998ExceptionalModulesTreeModules} that given a certain basis $\B$ of a representation $M$ of $Q$ (that is, a collection of basis elements from bases for vector spaces at all vertices), we say that the \emph{coefficient quiver $\Gamma(M,\B)$ of $M$ with respect to $\B$}\index{quiver!coefficient} is the quiver with vertex set $\B$ and having an arrow $b \xrightarrow{\alpha} b'$ provided the entry corresponding to $b$ and $b'$ in the matrix corresponding to $M(\alpha)$ with respect to the chosen basis $\B$ is non-zero.
Instead of labelling the arrows by the elements (of a $k$-basis) of $kQ$, we may also colour them accordingly and use the label to keep track of the corresponding coefficient of the structure matrix.

Further, we will call an indecomposable $k$-representation $M$ of $Q$ a \emph{tree module} provided there exists a basis $\B$ of $M$ such that the coefficient quiver $\Gamma(M, \B)$ is a tree.
Similarly, we will call $M$ \emph{planar} provided there exists a basis $\B$ such that $\Gamma(M,\B)$ is a planar graph.

Given a homomorphism $f \colon M \to N$, we will also consider the \emph{mapping quiver of~$f$}, denoted $\Gamma(f, \B_M, \B_N)$. This is a bipartite graph with vertex set $\B_M \cup \B_N$ and having an arrow $b \xrightarrow{f} b'$ provided the entry corresponding to $b$ and $b'$ in the structure matrix of $f$ wrt.\ to the given bases is non-zero.
Further, we say that $f$ is a \emph{tree map with respect to bases $\B_M$ and $\B_N$} provided $\Gamma(f,\B_M,\B_N)$ is a tree.

Let us now also recall the following property of coefficient quivers.

\begin{lemma}[{\cite[Property~1]{Ringel1998ExceptionalModulesTreeModules}}]
If $M$ is indecomposable and $\B$ is a basis of $M$, then $\Gamma(M,\B)$ is connected.

If $M$ is decomposable, then there exists a basis $\B$ of $M$ such that $\Gamma(M,\B)$ is not connected.
\end{lemma}

\subsection{Graph-theoretic background}

Recall that a graph $G$ is given by its set of vertices $V$ and a set of edges $E$ containing ordered pairs $(u,v) \in V^2$, describing an edge starting at $u$ and ending at $v$.

\textcite{EdwardsMcDiarmid1994NewUpperBoundHarmoniousColorings} have introduced a notion for classes of graphs which is similar to hyperfiniteness for families of modules.

\begin{mydef}[\cites{EdwardsMcDiarmid1994NewUpperBoundHarmoniousColorings}]
We say that a class $\Gamma$ of graphs is \defined{fragmentable} provided for any non-negative real number~$\eps$ there are positive integers $n_0$ and $C(\eps)$ such that for each graph $G \in \Gamma$ with $n\geq n_0$ non-isolated vertices, there is a set $X \subseteq V(G)$ of vertices such that
\begin{enumerate}
\item $\card{X} \leq \eps n$, and
\item each component of $G \setminus X$ has at most $C(\eps)$ vertices.
\end{enumerate}
\end{mydef}\index{fragmentability}


\begin{proposition} \label{prop:FragmentableCoeffQuiverBoundIndegreeHF} \label{prop:PlanarCoeffQuiverBoundIndegreeHF}
Let $d, \ell \in \N$.
Let $A$ be a finite dimensional $k$-algebra. 
Let $\M$ be a class of indecomposable $A$-modules, such that the class of underlying graphs of their coefficient quivers $\Gamma$ is fragmentable, their indegree is bounded by $d$ and the path length is bounded by $\ell$.
Then $\M$ is hyperfinite.
\begin{proof}
Let $\eps > 0$. Set $\tilde{\eps} = \eps \left(\sum_{i=0}^{\ell} d^i\right)^{-1}$ and pick $L_\eps \coloneqq c(\tilde{\eps})$, the constant from the definition of fragmentability.
Let $M \in \M$.
By the definition of fragmentability, there exists a set $S$ of vertices of the coefficient quiver of $M$ of cardinality at most $\tilde{\eps} \dim M$, such that the underlying graph $G$ splits into components of size at most $L_\eps$ if we remove the $\card{S}$ vertices. Now, if all the vertices in $S$ are sources, this subgraph describes a submodule.
As this might not hold, for each $s \in S$, we also remove all the vertices that map to $s$, of which there are at most $\sum_{i=1}^{\ell}d^i$. In this way, we arrive at a submodule $N$ of dimension at least $\dim M - \tilde{\eps} \sum_{i=0}^{\ell} d^i \dim M$. This finishes the proof.
\end{proof}
\end{proposition}

\begin{corollary} \label{cor:PlanarCoeffQuiverBoundIndegreeHF}
Let $d, \ell \in \N_0$.
Let $\M$ be the family of all the modules that have a planar coefficient quiver for which the indegree is bounded by $d$ and the path length bounded by $\ell$.
Then $\M$ is hyperfinite.
\begin{proof}
Planar graphs have genus $\gamma = 0$ and their subgraphs are also planar. So we can apply \cite[Corollary~3.7 (ii)]{EdwardsMcDiarmid1994NewUpperBoundHarmoniousColorings}---using a planar separator theorem. 
This shows that the family $\M$ is fragmentable.
Now apply Proposition~\ref{prop:FragmentableCoeffQuiverBoundIndegreeHF} and make use of the fact that subgraphs still satisfy the given bounds.
\end{proof}
\end{corollary}

\section{String algebras}

Given a quiver $Q = (Q_0,Q_1,s,t)$, we can consider an ideal $P$ of paths in $Q$ of length at least two. Then we can consider the quotient algebra $kQ/(P)$.

\begin{mydef}[{\cite{ButlerRingel1987AuslanderReitenSequencesFewMiddleTermsApplicationsStringAlgebras}}]
Let $A = kQ/P$. Then $A$ is a \defined{string algebra} provided
\begin{enumerate}
\item[(S1)] at most two arrows start in any vertex, and at most two arrows terminate in any vertex, i.e.\ $\card*{\set{a \in Q_1 \suchthat s(a) = v}}, \card*{\left\{a \in Q_1 \suchthat t(a) = v\right\}} \leq 2$ for all $v \in Q_0$,
\item[(S2)] for any $b \in Q_1$, there is at most one $a \in Q_1$ with $ba \notin P$, and at most one $c \in Q_1$ with $cb \notin P$,
\item[(S3)] for any $b \in Q_1$, there is some bound $n(b)$ such that any path $b_1 \dots b_{n(b)}$ with $b_1 = b$ contains a subpath in $P$ and there is some bound $n'(b)$ such that any path $b_{n'(b)} \dots b_1$ with $b_1 = b$ contains a subpath in $P$.
\end{enumerate}
\end{mydef}

\begin{example*}
Gentle algebras are a subclass of string algebras.
More concretely, the Kronecker algebra $k\Theta(2)$, i.e.\ the path algebra of the Kronecker quiver $1 \doublerightarrow 2$ is a string algebra. In the following, we shall denote its two arrows by $a$ and $b$.
\end{example*}

A \emph{letter} for $A$ is an arrow $a \in Q_1$ (a \emph{direct letter}) or its formal inverse $a^{-}$ (\emph{inverse letter}).
Start and terminus for letters are defined in the obvious way.
A \emph{word} $w$ is a formal sequence $w_1 \dots w_n$ of letters with $s(w_i) = t(w_{i+1})$, and we also define the inverse word $w^{-}$ in the apparent way.
A word of length $n$ is said to be a \emph{string} provided $w_i^{-} \neq w_{i+1}$ and for no subword $u = w_{i} \dots w_{i+t}$ we have $u,u^{-} \in P$.
A string $w$ is said to be a \emph{band} provided all powers $w^\ell = w \cdots w$ are also strings and $w$ is not itself the (proper) power of a shorter string.
Denote the set of all strings by $St$ and that of all bands by $Ba$.

\begin{example*}
For $\ell \geq 0$, strings for the Kronecker algebra include $a (b^{-}a)^\ell b^{-}$ and $b^{-} (ab^{-})^\ell a$. A band is given by $ab^{-}$.
\end{example*}

We use these words to define indecomposable $A$-modules, the \emph{string} and \emph{band modules}:

\begin{itemize}
\item[(string module)] For a string $s = s_1 \dots s_n$ we define a module $M(s)$ in the following way: As a vector space, $M(s) \isom \oplus_{j=0}^{n} V_j$, where $V_j = k$. The (non-trivial) action of $a \in Q_1$ is given by \begin{gather*}a(v_j) = \begin{cases}
 v_{j-1}, & j\geq 1, \, s_{j} = a,\\
 v_{j+1}, & j\leq n-1, \, s_{j+1} = a^{-},\end{cases}\end{gather*} where each $V_j = \spann{v_j}$.
 
\item[(band module)] Similarly, for a band $b = b_1 \dots b_n$ and an indecomposable $k[T,T^{-1}]$-module $(V,\varphi)$ of dimension $m$, we define a module $M(b,\varphi)$ in the following way: As a vector space, $M(b,\varphi) \isom \oplus_{0}^{n-1} V_j$, where $V_j = V$. The (non-trivial) action of $a \in Q_1$ is given by \begin{gather*}a(v_{i,j}) = \begin{cases}
 v_{i,j-1}, & 2 \leq j \leq n, b_j = a, \\
 v_{i,j+1}, & 1 \leq j \leq n-1, b_{j+1} = a^{-}, \\
 \varphi(v_{i,0}), & j=1, b_1 = a, \\
 \varphi^{-1}(v_{i,1}), & j=0, b_1 = a^{-}, \end{cases}\end{gather*} where each $V_j = \spann{v_{i,j}}_{i=1,\dots,m}$ and the values of the second index $j$ are treated modulo $n$.
\end{itemize}
 
Note that in the above cases, the action of the (remaining) idempotents $\epsilon_i$ of $A$ is the identity action on the basis elements corresponding to the vertex~$i$.

\begin{example}
The string module $M = M(ab^{-}ab^{-})$ has coefficient quiver \[v_0 \xleftarrow{a} v_1 \xrightarrow{b} v_2 \xleftarrow{a} v_3 \xrightarrow{b} v_4,\]
where $M_1 = \spann{v_1,v_3}$ and $M_2 = \spann{v_0,v_2,v_4}$, when we consider $M$ as a representation of the Kronecker quiver $\Theta(2) = 1 \doublerightarrow 2$.

For $\lambda \neq 0$, the band module $M\left(ab^{-},J_3(\lambda)\right)$, has coefficient quiver

\begin{center}
\begin{tikzcd}
{v_{3,0}}                                                         &  & {v_{2,0}}                                                                          &  & {v_{1,0}}                                                                          \\
{v_{3,1}} \arrow[u, "a", shift left] \arrow[u, "b"', shift right] &  & {v_{2,1}} \arrow[u, "a", shift left] \arrow[llu, "a"] \arrow[u, "b"', shift right] &  & {v_{1,1}}. \arrow[u, "a", shift left] \arrow[llu, "a"] \arrow[u, "b"', shift right]
\end{tikzcd}\end{center}

\noindent As a Kronecker representation, we have $M_1 = \spann{v_{i,1}}_i$ and $M_2 = \spann{v_{i,0}}_i$
Here, $J_n(\lambda)$ denotes the Jordan block of eigenvalue $\lambda$ for vector spaces $V_j = k^n$, given with respect to the basis $\{v_{1,j},\dots,v_{n,j}\}$.
\end{example}

Indecomposable modules for string algebras are classified by string and band modules. For this classification, let $\underline{\Phi}$ be a set of representatives of indecomposable $k[T,T^{-1}]$-modules.
\begin{theorem*}[{\cites[Proposition~2.3]{WaldWaschbuesch1985TameBiserialAlgebras}[Section~3]{ButlerRingel1987AuslanderReitenSequencesFewMiddleTermsApplicationsStringAlgebras}}]
Given a string algebra $A$, the string modules $M(w)$ for $w \in St$, and the band modules $M(w,\varphi)$ with $w \in Ba$, $\varphi \in \underline{\Phi}$ provide a complete list of indecomposable, finite dimensional $A$-modules.
\end{theorem*}

\begin{proposition} \label{prop:SAStringModulesHyperfinite}
The family of (indecomposable) string modules $\set{M(s) \colon w \in St}$ is hyperfinite.
\begin{proof}
We observe (as in \cite[Remark~3]{Ringel1998ExceptionalModulesTreeModules}) that the string modules $M(s)$ are indeed tree modules: by the construction of $M(s)$ for $s = s_1 \dots s_n$, $M(s)$ is of dimension $n+1$ and for each $s_i$ there is exactly one edge in the coefficient quiver. Hence there are $n$ edges. But the coefficient quiver is connected since the module is indecomposable. 
Indeed, string modules are also tree modules in the ``stronger'' sense of \cite[Section~1]{CrawleyBoevey1989MapsBetweenRepresentationsZeroRelationAlgebras}, i.e.\ at most one arrow of each label/colour starts at each vertex. 
To see this, acknowledge that in order to have an edge $v_j \edge v_k$ in $\Gamma(M(s),\B)$, we must have $\abs{j-k} = 1$, and each of these edges is determined by $s_{\max(j,k)}$, hence there can be at most two edges starting or terminating at any given $v_j$, and if $v_j$ is a sink our source, they cannot have the same label since consecutive letters cannot be each other's inverse.
This implies that the indegree is globally bounded. 
Hence, we may apply Proposition~\ref{prop:FragmentableCoeffQuiverBoundIndegreeHF} to see that the indecomposable string modules form a hyperfinite family.
\end{proof}
\end{proposition}

We shall now use this result for string modules to give a proof different from \cite{Elek2017InfiniteDimensionalRepresentationsAmenabilty} for the amenability of string algebras by reducing the study of band modules to string modules.

\begin{theorem}[{\cite[Proposition~10.1]{Elek2017InfiniteDimensionalRepresentationsAmenabilty}}]
Let $A$ be a finite dimensional $k$-algebra. If $A$ is a string algebra [special biserial algebra], then $A$ is of amenable representation type.
\begin{proof}
According to the classification theorem, every indecomposable module of a string algebra is either a string or a band module.
By Proposition~\ref{prop:SAStringModulesHyperfinite}, the family of indecomposable string modules is hyperfinite. It therefore remains to consider the band modules.
By the proof of \cite[Proposition~18.1]{Ringel2011MinimalRepresentationInfiniteAlgebrasSpecialBiserial}, each indecomposable band module has a maximal submodule $N$ which is a string module.
Indeed, the submodule constructed there---as a vector space---is (in our notation) $N = \oplus_{j=0}^{n-2}V_j \oplus U$, where $U$ is a direct complement in $V_{n-1}$ of a one-dimensional subspace. 
Hence, each band module has a submodule of codimension $1$, which is a string module. We may thus apply Proposition~\ref{prop:ExtendingHFfromSubmodulesOfBoundedCodimension}.

To also include special biserial algebras, note that in addition to tree and band modules, there are no more (finite dimensional) indecomposable modules except some projective-injective modules, but of those there can only be finitely many.
\end{proof}
\end{theorem}

\section{Clannish algebras}

To prove amenability of clannish algebras, we want to mimic the proof for string algebras.
Here, we need to deal with the fact that the classification of the indecomposable modules is slightly more complicated, but was achieved by \cites{CrawleyBoevey1989ClansGelfandProblem}{Bondarenko1991RepresentationsBundlesSemichainedSetsTheirApplications}{Deng2000ProblemNazarovaRoiter}{Hansper2022PhDThesis}.

We will recall the relevant definitions and the results from \cites{CrawleyBoevey1989ClansGelfandProblem}{Geiss1999MapsBetweenRepresentationsOfClans} and \cite{Hansper2022PhDThesis}.

Let $Q=(Q_0,Q_1,s,t)$ be a quiver, $\Sp \subseteq Q_1$ a subset of the loops. We say that $\eps \in \Sp$ is a \emph{special arrow} (or \emph{loop}), while we call the remaining arrows $a \in Q_1 \setminus \Sp$ \emph{ordinary arrows}, denoting their set by $\Ord$. 
Let $R^\Sp \coloneqq \set*{\eps^2 - \eps \suchthat \eps \in \Sp}$ be the idempotent relations for the special arrows and let $R$ be a set of monomial relations on $Q$. 

\begin{mydef}[{\cites[Section~2.5]{CrawleyBoevey1989ClansGelfandProblem}[5.6 Definition]{Geiss1999MapsBetweenRepresentationsOfClans}}]
Let $\Lambda = kQ/(R \cup R^\Sp)$. We say that $\Lambda$ is \defined{clannish} provided the following conditions hold: 
\begin{enumerate}
\item[(C0)] no relation in $R$ starts or ends in a special loop or involves the square of one,
\item[(C1)] at most two arrows start at any vertex, and at most two arrows terminate at any vertex, i.e.\ $\card*{\set{a \in Q_1 \suchthat s(a) = v}}, \card*{\left\{a \in Q_1 \suchthat t(a) = v\right\}} \leq 2$ for all $v \in Q_0$,
\item[(C2)] for any $b \in \Ord$, there is at most one $a \in Q_1$ with $ba \notin R$ and at most one $c \in Q_1$ with $cb \notin R$.
\end{enumerate}
\end{mydef}

\begin{examples*}
String algebras and skewed-gentle algebras are clannish. 

The quiver $Q = \begin{tikzcd}& 1 \arrow["\eps", loop right, distance=2em] \arrow["a", loop left, distance=2em] \end{tikzcd}$ along with $\Sp = \{\eps\}$ and $R = \left\{a^2 = 0\right\}$ gives the clannish algebra $k\spann{a,\eps}/(a^2,\eps^2-\eps)$. It is not finite dimensional over $k$.

As in \cite[Example~2.14]{Hansper2022PhDThesis}, consider the quiver \[Q = \begin{tikzcd}[row sep=tiny, column sep=small]
\arrow["a", dr] \arrow["\eps", loop left]1 & & 5 \arrow["e", dd] \arrow["\kappa", loop right]\\
& 2 \arrow["c" ,ur] \arrow["d", dr]\\
\arrow["\eta", loop left] 3 \arrow["b", ur] & & 4
\end{tikzcd}.\]
Then $Q$ along with $\Sp = \set{\eps, \eta, \kappa}$ and $R = \set{ca, db, ec}$ gives a clannish algebra.
\end{examples*}


A \emph{letter} for $\Lambda$ is an arrow $a \in \Ord$ (a \emph{direct letter}), its formal inverse $a^{-}$ (\emph{inverse letter}), or a symbol $\eps^{*}$ for any $\eps \in \Sp$ (\emph{special letter} and its own inverse).
Source and target of letters are defined in the obvious way, as are formal inverses of (sequences of) letters.
A \emph{word} $w$ in $\Lambda$ is a (finite) formal sequence $w_1 w_2 \dots w_n$ of letters with $s(w_i) = t(w_{i+1})$ and $w_i^{-} \neq w_{i+1}$ such that $w$ contains no subword $u$ with $u, u^{-} \in R$ and for $\eps \in \Sp$, $\eps^* \eps^*$ is not a subsequence of $w$. Words can be composed if their concatenation itself is a word while the inverse of $w$ (as before) is $w_n^{-} \dots w_1^{-}$. 
We say that a word $w$ is \emph{coadmissible} provided it cannot be pre- or postcomposed with a special letter to give a longer word.
A \emph{rotation} of a word $w_1 \dots w_n$ is a word of the form $w_{i+1} \dots w_{n} w_{1} \dots w_{i}$ for some $i \in \set{1,\dots,n}$.

Let $St$  be a set of representatives of words $w$ which are coadmissible, under the equivalence relation which identifies $w$ with $w^{-}$. We call the elements of $St$ \emph{strings} and say that a string $w$ is \emph{symmetric} if $w = w^{-}$, (and then up to equivalence $w = zf^*z^-$ for some word $z=z_1 \dots z_n$ and special loop $f$) and otherwise \emph{asymmetric}.
Let $Ba$ be a set of representatives of nontrivial words $w$ such that $ww$ is defined and such that $w$ is not a power, under the equivalence relation which identifies all rotations of $w$ and $w^{-}$. We call the elements of $Ba$ \emph{bands} and say that a band $w$ is \emph{symmetric} if $w$ is equal to some rotation of $w^-$ (and then up to equivalence $w=f^{*}z^{-}g^{*}z$ for special loops $f,g$ and some word $z$), and otherwise \emph{asymmetric}.

\begin{example*}
Let $A=k\spann{a,\eps}/(a,\eps^2-\eps)$.
$\eps^{*}a\eps^{*}$ is a word in $A$ with inverse $\eps^{*}a^{-}\eps^{*}$.
This is an asymmetric string. Elements of $Ba$ include the symmetric band $\eps^{*}a^{-}\eps^{*}a$ and the asymmetric band $a\eps^{*}$.

On the other hand, for $A = kQ/(ca,db,ec,\eps^2-\eps,\eta^2-\eta,\kappa^2-\kappa)$ as above, the word $da\eps^*a^-d^-$ is a symmetric string with $z = da$ 
and the word $\eps^*a^-b\eta^*b^-a$ is a symmetric band with $z = b^-a$.
\end{example*}

Let us recall how we construct indecomposable (clannish) \emph{string} and \emph{band modules} for these algebras.

First, we have to use some additional information: We need to associate a \emph{direction} $\dir \colon \set{1,\dots,l} \to \set{\pm}$ to letters of a given word $w_1 \dots w_l$ of length $l$. For direct ($+$) and inverse ($-$) letters, this is straightforward. For special letters, the direction is determined by a partial order on the set of words, induced by a linear order on the letters. 
Roughly speaking, we have that $\dir(w_i) = +$ iff the inverse of the subword to the left of $w_i$ is---in this partial order---less than the subword to the right, taking into consideration the symmetry properties of the word (see, e.g., \cite[Chapter~3]{Hansper2022PhDThesis}).

Now, in order to construct (clannish) string and band modules, we define for each $w \in St \cup Ba$ an algebra $A_w$ and a functor $S_w \colon \mods A_w \to \mods A$.
For our purpose, it suffices to describe these functors on $\ind A_w$.

\begin{itemize}
\item[(asymmetric string module)] Given an asymmetric string $w = w_1 \dots w_n$, we put \mbox{$A_w = k$} and have $S_w(k) \isom \oplus_{j=0}^n V_j$ as vector spaces, where $V_j = k$. The (non-trivial) actions of $a \in \Ord$ and $\eps \in \Sp$ are given by \begin{gather*}a(v_j) = \begin{cases}
 v_{j-1}, & 
 w_{j\phantom{+1}} = a,\\
 v_{j+1}, & 
 w_{j+1} = a^{-},\end{cases} \quad \eps(v_j) = \begin{cases}
  v_{j-1}, & w_{j\phantom{+1}} = \eps^*, \dir(w_{j\phantom{+1}}) = +,\\
  v_{j+1}, & w_{j+1} = \eps^*, \dir(w_{j+1}) = -,\\
  v_{j}, & w_{j\phantom{+1}} = \eps^*, \dir(w_{j\phantom{+1}}) = -,\\
  v_{j}, & w_{j+1} = \eps^*, \dir(w_{j+1}) = +,\end{cases}\end{gather*} where each $V_j = \spann{v_j}$.
\item[(symmetric string module)] For a symmetric string $w = zf^{*}z^{-}$ with $z = z_1 \dots z_n$, we put $A_w = k[T]/(T^2-T)$ and have $S_w(V) \isom \oplus_{j=0}^{n} V_j$ as vector spaces, where $V_j = V = k$, either with $T$ acting on $V$ as the identity or the zero map. The (non-trivial) actions of $a \in \Ord$ and $\eps \in \Sp$ are given by \begin{gather*}a(v_j) = \begin{cases}
 v_{j-1}, & 
 z_{j\phantom{+1}} = a,\\
 v_{j+1}, & 
 z_{j+1} = a^{-},\end{cases} \quad
\eps(v_j) = \begin{cases}
   v_{j-1}, & z_{j\phantom{+1}} = \eps^*, \dir(z_{j\phantom{+1}}) = +,\\
   v_{j+1}, & z_{j+1} = \eps^*, \dir(z_{j+1}) = -,\\
   v_{j}, & z_{j\phantom{+1}} = \eps^*, \dir(z_{j\phantom{+1}}) = -,\\
   v_{j}, & z_{j+1} = \eps^*, \dir(z_{j+1}) = +,\\\hline
   T \cdot v_j, & j=n, f^* = \eps^*,\end{cases}
   \end{gather*}
    where each $V_j = \spann{v_j}$.
\item[(asymmetric band module)] For an asymmetric band $w = w_1 \dots w_n$, we put $A_w = k[T,T^{-1}]$ and have $S_w(V) = \oplus_{j=0}^{n-1} V_j$ as vector spaces, where $V_j = (V,\phi)$ is an indecomposable $A_w$-module of dimension $m$. The (non-trivial) actions of $a \in \Ord$ and $\eps \in \Sp$ are given by \begin{gather*}a(v_{i,j}) = \begin{cases}
 v_{i,j-1}, & j \neq 1, w_{j\phantom{+1}} = a, \\
 v_{i,j+1}, & j \neq 0, w_{j+1} = a^{-}, \\\hline
 \phi(v_{i,0}), & j=1, w_1 = a, \\
 \phi^{-1}(v_{i,1}), &j=0, w_1 = a^{-}, \end{cases} \quad
\eps(v_{i,j}) = \begin{cases} 
  v_{i,j-1}, & j \neq 1, w_{j\phantom{+1}} = \eps^*, \dir(w_{j\phantom{+1}}) = +,\\
  v_{i,j+1}, & j \neq 0, w_{j+1} = \eps^{*}, \dir(w_{j+1}) = -,\\
  v_{i,j}, & j \neq 1, w_{j\phantom{+1}} = \eps^*, \dir(w_{j\phantom{+1}}) = -,\\
  v_{i,j}, & j \neq 0, w_{j+1} = \eps^{*}, \dir(w_{j+1}) = +,\\\hline
  \phi(v_{i,0}), & j=1, w_{1\phantom{+1}} = \eps^*, \dir(w_{1\phantom{+1}}) = +,\\
  v_{i,0}, & j=0, w_{1\phantom{+1}} = \eps^*, \dir(w_{1\phantom{+1}}) = +\\
  \phi^{-1}(v_{i,1}), & j=0, w_{1\phantom{+1}} = \eps^*, \dir(w_{1\phantom{+1}}) = -,\\
  v_{i,1}, & j=1, w_{1\phantom{+1}} = \eps^*, \dir(w_{1\phantom{+1}}) = -\end{cases}\end{gather*} where each $V_j = \spann{v_{i,j}}_{i=1,\dots,m}$.

\item[(symmetric band module)] For a symmetric band $w = f^{*}z^{-}g^{*}z$ with $z = z_1 \dots z_n$, we put $A_w = k\langle x,y\rangle/(x^2-x,y^2-y)$ and have $S_w(V) = \oplus_{j=0}^{n} V_j$ as vector spaces, where $V_j = (V,\phi,\psi)$ is an indecomposable $A_w$-module of dimension $m$. The (non-trivial) actions of $a \in \Ord$ and $\eps \in \Sp$ are given by \begin{gather*}a(v_{i,j}) = \begin{cases}
 v_{i,j-1}, & 
 z_{j\phantom{+1}} = a, \\
 v_{i,j+1}, & 
 z_{j+1} = a^{-}, \end{cases} \quad
\eps(v_{i,j}) = \begin{cases}
  \psi(v_{i,0}), & j=0, g^{*} = \eps^*, \\\hline
  v_{i,j-1}, & z_{j\phantom{+1}} =\eps^*, \dir(z_{j\phantom{+1}}) = +,\\
  v_{i,j+1}, & z_{j+1} = \eps^*, \dir(z_{j+1}) = -,\\
  v_{i,j}, & z_{j\phantom{+1}} = \eps^*, \dir(z_{j\phantom{+1}}) = -,\\
  v_{i,j}, & z_{j+1} = \eps^*, \dir(z_{j+1}) = +,\\\hline
  \phi(v_{i,n}), & j=n, f^* = \eps^*, \end{cases}\end{gather*} where each $V_j = \spann{v_{i,j}}_{i=1,\dots,m}$.
\end{itemize}
Note that in the above cases, the action of the (remaining) idempotents $\epsilon_i$ of $A$ is the identity action on the basis element corresponding to the vertex $i$.

\begin{example}
Consider $A=kQ/(ca,db,ec,\eps^2-\eps,\eta^2-\eta,\kappa^2-\kappa)$ as before.
Let $V = k$ be the $k[T]/(T^2-T)$-module with $T$ acting as the identity. Then the symmetric string module $S_{\eps^{*}a^{-}b\eta^*b^-a\eps^*}(V)$ has coefficient quiver
\[\begin{tikzcd}
v_0 \arrow["\eps", r] & v_1 \arrow["\eps", loop, distance=2em, in=125, out=55] \arrow["a", r] & v_2 & v_3 \arrow["b"', l] \arrow["\eta", loop right, distance=2em]
\end{tikzcd}.\]

Now consider $k\spann{a,\eps}/(a,\eps^2-\eps)$.
Let $\left(k^3,\begin{psmallmatrix}0&0&0\\1&1&0\\1&0&1\end{psmallmatrix},\begin{psmallmatrix}1&0&0\\0&1&0\\0&0&0\end{psmallmatrix}\right)$ be an indecomposable $k\spann{x,y}/(x^2-x,y^2-y)$-module. The symmetric band module $S_{\eps^{*}a^{-}\eps^{*}a}(V)$ has coefficient quiver wrt.\ to the basis $\set{v_{i,j}}_{i=1,\dots,3; j=0,1}$ 
\[\begin{tikzcd}
{v_{3,0}} & {v_{2,0}} \arrow["\eps", loop, distance=2em, in=125, out=55]            & {v_{1,0}} \arrow["\eps", loop, distance=2em, in=125, out=55] \\
          & {v_{2,1}} \arrow["a", u] \arrow["\eps", loop left, distance=2em] &   \\
{v_{3,1}} \arrow["a", uu] \arrow["\eps", loop left, distance=2em] &      & {v_{1,1}.} \arrow["a", uu] \arrow[lu, "\eps"] \arrow[ll, "\eps"]
\end{tikzcd} 
\]
\end{example}

Let us now recall the classification of the indecomposable modules of a clannish algebra.

\begin{theorem*}[{\cites[3.8 Main Theorem]{CrawleyBoevey1989ClansGelfandProblem}[Theorem~1]{Geiss1999MapsBetweenRepresentationsOfClans}[Main Theorem]{Hansper2022PhDThesis}}]
If $A = kQ/(R\cup R^\Sp)$ is a clannish algebra, the representations $S_w(V)$ with $w$ running through $St \cup Ba$ and $V$ running through $\ind A_w$, provide a complete set of non-isomorphic indecomposable representations of $A$.
\end{theorem*}

We can now consider the different types of modules with respect to hyperfiniteness.

\begin{proposition} \label{prop:ClanAlgStModulesPlanar}
Let $A$ be a clannish algebra.
If $M$ is an indecomposable clannish string module,
there is a $k$-basis $\B$ of $M$ such that the coefficient quiver $\Gamma(M,\B)$ is a planar graph and has indegree bounded by $3$.
\begin{proof}
For the string module case, similar to string modules for string algebras, the action of ordinary arrows effects edges $v_j \edge v_{j\pm1}$, but at most two for each~$j$.
Special loops may additionally add up to two loops $v_j \rightloop$, breaking neither planarity nor indegree bound. 
We may thus simply choose $\B = \set{v_0,\dots,v_n}$ as in the construction.
\end{proof}
\end{proposition}

\begin{lemma} \label{lemma:PlanarCoeffQuiverInvMap}
Let $\Lambda = k[T,T^{-1}]$ and let $M = (V,\phi)$ be an indecomposable finite dimensional $\Lambda$-module.
Then $V$ has a subspace $U$ of codimension $1$ with a planar mapping quiver $\Gamma(\phi_{|U},\B_U,\B_V)$. 
\begin{proof}
A finite dimensional $\Lambda$-module is given by a finite dimensional vector space $V$ along with an endomorphism $\phi$. Since $\Lambda$-isomorphisms translate to replacing $\phi$ by a similar matrix, we may assume that $\phi$ is given in rational canonical form. Since $(V,\phi)$ is also assumed to be indecomposable, $\phi$ is just the companion matrix of a monic polynomial. Choose a basis $\set{v_1,\dots,v_m}$ of $V$ corresponding to $\phi$, where $v_m$ is associated to the column of $\phi$ containing the coefficients of the polynomial. Now put $U = \spann{u_1,\dots,u_{m-1}}$, where $u_i = v_i$. Then $\Gamma(\phi_{|U},\B_U,\B_V)$ is planar, indeed a forest containing just the trees $u_{i} \to v_{i+1}$ for $1\leq i \leq m-1$.
\end{proof}
\end{lemma}

\begin{proposition} \label{prop:ClanAlgAsymBaModulesPlanar}
Let $A$ be a finite dimensional clannish algebra.
If $M$ is an asymmetric band module,
there is a submodule $N\subset M$ of codimension $1$ such that $N$ has a $k$-basis $\B$ such that the coefficient quiver $\Gamma(N,\B)$ is a planar graph and has indegree bounded by $4$.
\begin{proof}

\begin{figure}
\begin{tikzcd}[row sep=1em,column sep=1.5em,/tikz/execute at end picture={
    \node (large) [rectangle, draw, color=red, fit=(A1) (A2)] {};
  }]
  & {v_{1,2}} \arrow[ldddd, no head, bend right=10] & |[alias=A1]| {u_1} \arrow[rrd] \arrow[l] &  & {v_{1,0}} \arrow[r, no head] & {v_{1,n-1}} \arrow[rdddd, no head, bend left=10] & \\
  & \vdots \arrow[lddd, no head, bend right=8]      & \vdots \arrow[l] \arrow[rrd] & \vdots & \vdots \arrow[r, no head] & \vdots \arrow[rddd, no head, bend left=8] & \\
  & {v_{m-1,2}} \arrow[ldd, no head, bend right=10] & u_{m-1} \arrow[l] \arrow[rrd] & & {v_{m-1,0}} \arrow[r, no head] & {v_{m-1,n-1}} \arrow[rdd, no head, bend left=10] & \\
  & {v_{m,2}} \arrow[ld, no head, bend right=5]     &  &  & |[alias=A2]| {v_{m,0}} \arrow[r, no head] & {v_{m,n-1}} \arrow[rd, no head, bend left=5] & \\
\B_{V_3} \arrow[rdd, no head, double equal sign distance, bend right] &  &  &  &  &  & \B_{V_{n-2}} \\
  &  &  &  &  &  &  \\
  & \B_{V_4} \arrow[rr, no head, double equal sign distance] &  & \dots \arrow[rr, no head, double equal sign distance] &  & \B_{V_{n-3}} \arrow[ruu, no head, double equal sign distance, bend right] & 
\end{tikzcd}
\caption{Sketch of the coefficient quiver of a general asymmetric band module, with a box showing the mapping quiver $\Gamma(\phi_{|U},\B_U,\B_V)$, omitting the possible loops due to special letters.} \label{fig:CoeffQuiverAsBa}
\end{figure}
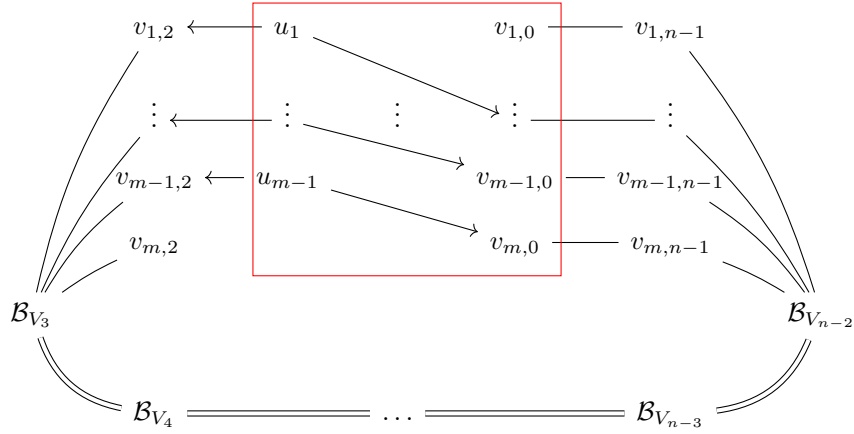

For asymmetric band modules $S_w(V,\phi)$, letters other than $w_1$ effect at most two edges $v_{i,j} \edge v_{i,j\pm1}$ and at most two loops $v_{i,j} \rightloop$.
By rotating $w$, we may assume that $\dir(w_1) = +$ and $\dir(w_2) = -$, i.e.\ wrt.\ to the actions of the algebra generators, $V_1$ is a source. Otherwise, all letters of $w$ would have the same direction, implying that $w$ resp.\ $w^{-}$ and their powers give paths in the quiver, contradicting the finite dimensionality of $A$.
It follows that the basis elements $\set{v_{i,1}}$ are mapped to $\set{v_{i,0}}$ according to the entries of the transformation matrix of $\phi$ wrt.\ to the basis $\B_V = \set{v_1,\dots,v_m}$ of $V_1$.
Construct a submodule by replacing the vector space $V_1$ by the space $U$ of Lemma~\ref{lemma:PlanarCoeffQuiverInvMap}. The action of $\phi$ wrt.\ to the basis $B_U = \set{u_1,\dots,u_{m-1}}$ of this subspace is given by a planar mapping quiver.
Hence, if we arrange the bases of the other $V_j$s in a circle with parallel edges realising the planarity (with some loops added), the remaining ``gap'' between $U$ and $V_0$ can just be filled with $\Gamma(\phi_{|U},\B_U,\B_V)$, keeping the coefficient quiver wrt.\ $\B =\set{v_{1,0},\dots,v_{m,0},u_1,\dots,u_{m-1},v_{1,2},\dots,v_{m,n-1}}$ planar (see Figure~\ref{fig:CoeffQuiverAsBa}). Also note that the mapping quiver section increases the indegree of each vertex by at most one.
Note that the submodule constructed in this way has codimension one.
\end{proof}
\end{proposition}

\begin{lemma} \label{lemma:PlanarInclusionTwoIdempotentMaps}
Let $\Lambda = k\langle x,y\rangle/(x^2-x,y^2-y)$ and let $M = (V,\phi,\psi)$ be an indecomposable $\Lambda$-module.
Then $V$ has a subspace $U$ of codimension at most $1$ and there are bases $\B_x$ and $\B_y$ of $V$ that respect the idempotent decompositions of $V$ induced by $\phi$ and $\psi$ respectively, and the inclusion $U \embeds V$ is a tree map with respect to $\B_y'$ and $\B_x$, where $\B_y' \subseteq \B_y$ with $U = \spann*{\B_y'}$.
\begin{proof}
Since $\phi$ and $\psi$ are idempotent endomorphisms of $V$, they induce two vector space decompositions $V \isom \ker\phi \oplus \im\phi \isom \ker\psi \oplus \im\psi$. Hence we can associate to $(V,\phi,\psi)$ a representation of the four subspace quiver $\S(4)$ in which the subspaces are paired and each pair is complimentary, namely \[\begin{tikzcd}
                          & \ker\phi \arrow["A", d, hook] &                          \\
\ker\psi \arrow["D", r, hook] & V                         & \im\phi. \arrow["B", l, hook] \\
                          & \im\psi \arrow["C", u, hook]   &                         
\end{tikzcd}\]
By the classification of indecomposable $k\S(4)$-modules, in this way we only get regular representations, and they lie in all but one of the tubes (using the notation of \cite[Section~XIII.3]{SimsonSkowronski2007TubesConcealedAlgebrasEuclideanType}, we miss the tube of index~$1$).
What is more, we can assume that $A,B,C,D$ are block matrices containing identity matrices (enlarged by zero columns or rows), a Jordan block of eigenvalue zero or the Frobenius companion matrix of a power of a monic, irreducible polynomial (see \cite[Appendix]{MedinaZavadskij2004FourSubspaceProblemElementarySolution} for an explicit description of the matrices).

We know that either $A, B$ or $C, D$ are $2 \times 1$-block matrices with their only-non zero block an identity matrix. We may hence assume that $V = \ker\phi \oplus \im\phi$ and put $\B_x$ to be the union of the canonical bases of $\ker\phi$ and $\im\phi$.
Let $\B_y$ be the union of the canonical bases of $\ker\psi$ and $\im\psi$.
Now, the block matrix $(C|D)$ serves as the basis change matrix $\B_y \to \B_x$.
If $(C|D)$ contains a companion matrix block (i.e.\ the representation is in a tube of rank one), 
we remove from $\B_y$ the basis element corresponding to the last column of the companion matrix and call the resulting set $\B_y'$ and denote the modified matrix by $(\widehat{C|D})$.
Otherwise, we put $\B_y' = \B_y$.

Next, let $U \coloneqq \spann{\B_y'}$ and consider the map $U \embeds V$ and its mapping quiver with respect to $\B_y'$ and $\B_x$. Its edges are determined by the matrix $(C|D)$ or $(\widehat{C|D})$. But by the structure of $C$ and $D$, we know that there are $\dim U + \dim V - 1$ edges in a graph which has $\dim U + \dim V$ vertices. Moreover, the degree bound is two. It follows that the mapping quiver is a bipartite tree, indeed a line.
To be more precise, depending on the type of the representation (as in \cite{MedinaZavadskij2004FourSubspaceProblemElementarySolution}), the edges are described by
\begin{enumerate}
\item[(0)] $e_i \from f_i \to e_{n+i}$ for $1\leq i \leq n$ and $e_{i+1} \from f_{n+i} \to e_{n+i}$ for $1\leq i \leq n-1$,
\item[(I)] $e_i \from f_i \to e_{n+i}$ for $1\leq i \leq n$, $f_{n+1} \to e_{n+1}$ and $e_{i-1}\from f_{n+i} \to e_{n+1}$ for $2\leq i \leq n+1$,
\item[(II)] $e_i \from f_i \to e_{n+1+i}$ for $1 \leq i \leq n$, $f_{n+1} \to e_1$ and $e_{1+i}\from f_{n+1+i} \to e_{n+1+i}$ for $1 \leq i \leq n$,
\end{enumerate}
for suitable enumerations $\B_x = \set{e_i}$ and $B_y' = \set{f_i}$.

Note that realising the planarity, e.g.\ when putting the basis elements on two opposing sides of a rectangle, induces (up to reversing the direction) a pair of canonical (re-)orderings of $\B_x$ and $\B_y'$.
\end{proof}
\end{lemma}

\begin{lemma} \label{lemma:SymBaUnidirectedFiniteLength}
Let $A$ be a finite dimensional clannish algebra.
Then there is some $L > 0$ such that any symmetric band containing direct and inverse letters, of length at least $2L+2$, contains letters of positive and negative direction within the first $L+1$ letters.
\begin{proof}
Towards a contradiction and without loss of generality, we may assume that for all $N \in \N$ there is $L \geq N$ such that there exists a symmetric band $w = f^{*} z^{-} g^{*} z$ such that $z$ has length $L$ and $\dir(z_i) = +$ for $1\leq i \leq L$.
By the construction of the modules associated to bands, this means that there is $b \in S_w(V,\phi,\psi)$ such that
$ z_1 \circ \dots \circ z_L(b)\neq 0$.
Hence, the path associated to this composition of structure maps is also non-zero. It has length $L$. Since $z$ does not involve subwords contained in $R$ and no squares of special loops, the path does not equal a shorter path. It follows that there are infinitely many such paths, which are linearly independent. But $A$ was assumed to be finite dimensional. A contradiction.
\end{proof}
\end{lemma}

\begin{proposition} \label{prop:ClanAlgSymBaModulesPlanarSubmodule}
Let $A$ be a finite dimensional clannish algebra. Then there is some $L > 0$ such that all indecomposable (clannish) symmetric band modules $M$ have a submodule $N$ of codimension at most $L$ such that $N$ has a basis $\B$ such that the coefficient quiver $\Gamma(N,\B)$ is planar and has indegree bound by $5$.
\begin{proof}
Let $M = S_{w}(V,\phi,\psi)$ be an indecomposable symmetric band module with $w = f^{*}z^{-}g^{*}z$ for some special loops $f,g$ and a word $z = z_1 \dots z_n$. Since $A$ is finite dimensional, $z$ cannot be empty or only contain special letters, as otherwise the powers of (the special loop path corresponding to) $w$ would give infinitely many linearly independent elements. 

Use Lemma~\ref{lemma:PlanarInclusionTwoIdempotentMaps} to construct a mapping quiver for $(V,\phi,\psi)$. We may assume that the basis $\B_y'$ is associated to $\psi$ [to $\phi$]. 
Hence, let $i_0$ be the smallest index such that $\dir(z_{i_0})=-$ [$\dir(z_{n-i_0})=+$]. If this cannot be satisfied, $z$ contains no inverse [direct] letters. Then, we instead consider $g^{*}z'^{-}f^{*}z'$ with $z' = z^{-}$, which gives an isomorphic clannish module, while $z'$ now contains an inverse [a direct] letter. 
By the previous Lemma~\ref{lemma:SymBaUnidirectedFiniteLength}, we must have $i_0 \leq L$, where the constant is the one from the Lemma.

Now we construct $N$ from $M$ by replacing the vector spaces $V_0, \dots, V_{i_0}$ {\linebreak} [$V_{n-i_0}, \dots, V_{n}$] by the subspace $U = \spann{B_y'}$ given by the above lemma.
\begin{figure}
\begin{tikzcd}
\B_y' \arrow[r, no head, double equal sign distance] \arrow["\psi"', loop below, distance=2em, out=255, in=285] & \B_y' \arrow[r, no head, double equal sign distance] \arrow[dotted, loop above, distance=2em, in=110, out=70] & \dots \arrow[r, no head, double equal sign distance] & \B_y' \arrow[r, "\left(\widehat{C|D}\right)", rightarrow] \arrow[dotted, loop above, distance=2em, in=110, out=70] & \B_x \arrow[r, no head, double equal sign distance] \arrow[dotted, loop above, distance=2em, in=70, out=110] & \dots \arrow[r, no head, double equal sign distance] & \B_x \arrow[r, no head, double equal sign distance] \arrow[dotted, loop, distance=2em, in=70, out=110] & \B_x \arrow["\phi", loop below, distance=2em, in=255, out=285]\\
&&&&\text{at $i_0$} \arrow[u,>=Stealth]
\end{tikzcd}
\caption{Coefficient quiver of a planar submodule of an indecomposable symmetric band module} \label{fig:CoeffQuiverSubModSymBa}
\end{figure}
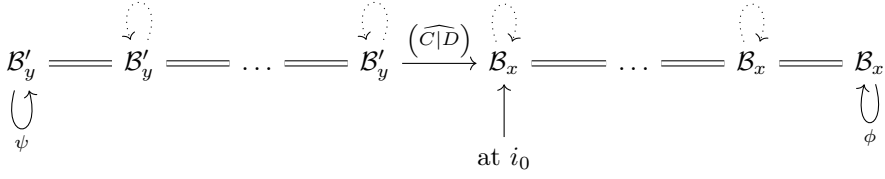
The coefficient quiver of $N$ has shape as depicted in Figure~\ref{fig:CoeffQuiverSubModSymBa}. Note that the edges between a basis of $U$ at $i_0-1$ [$n-i_0+1$] and of $V$ at the neighbouring vertex are given by the mapping quiver obtained previously. The reordering needed does not affect the planarity in general, as all other edges are of the form $v_{i,j} \edge v_{i,j\pm1}$ or are loops.
Moreover note that the bases $\B_x$ and $\B_y$ were chosen in such a way that the action of $\phi$ and $\psi$ is either a loop $v_{n,j} \rightloop$ respectively $v_{0,j} \rightloop$ or zero.
Hence the coefficient quiver is a planar graph with degree bounded by five, 
and the path length is bounded by $\dim_k A$.
We also note that $\dim_k M -\dim_k N \leq L$, since we only change the vector spaces at a certain number of vertices.
\end{proof}
\end{proposition}

\begin{theorem}
Let $A$ be a finite dimensional, clannish $k$-algebra. Then $A$ is of amenable representation type.
\begin{proof}
By the classification theorem, every indecomposable $A$-module is a symmetric or asymmetric string or band module.
By Proposition~\ref{prop:ClanAlgStModulesPlanar}, clannish string modules are planar with an upper bound on degree and path length. To their family, we can therefore directly apply Proposition~\ref{prop:PlanarCoeffQuiverBoundIndegreeHF}.
Also by Proposition~\ref{prop:ClanAlgAsymBaModulesPlanar}, asymmetric band modules have submodules of globally bounded codimension with planar coefficient quivers and a degree bound. Here, we therefore additionally apply Proposition~\ref{prop:ExtendingHFfromSubmodulesOfBoundedCodimension}.
Finally, for symmetric band modules, we apply Proposition~\ref{prop:ClanAlgSymBaModulesPlanarSubmodule} to see that they have submodules of globally bounded codimension with planar coefficient quivers and a degree bound. Again, their hyperfiniteness follows from an application of Propositions~\ref{prop:PlanarCoeffQuiverBoundIndegreeHF} and~\ref{prop:ExtendingHFfromSubmodulesOfBoundedCodimension}.
\end{proof}
\end{theorem}

{
 \printbibliography
}
\end{document}